\documentclass{article}
\usepackage{amsfonts,amssymb}
\usepackage{amsthm}
\usepackage[centertags]{amsmath}
\def\bf{\bfseries}
\def\it{\itshape}

\newcommand{\G}{{\cal G}}

\newcommand{\eps}{\varepsilon}

\newtheorem{theorem}{Theorem}
\newtheorem{itlemma}{Lemma}[section]
\newtheorem{itclaim}{Claim}[section]

\newtheorem{itproposition}[itlemma]{Proposition}
\newtheorem{itcorollary}[itlemma]{Corollary}
\newtheorem{itremark}[itlemma]{Remark}
\newtheorem{itdefinition}[itlemma]{Definition}
\newtheorem{itexample}[itlemma]{Example}
\newenvironment{claim}{\begin{itclaim}\rm}{\end{itclaim}}
\newenvironment{lemma}{\begin{itlemma}\rm}{\end{itlemma}}
\newenvironment{remark}{\begin{itremark}\rm}{\end{itremark}}
\newenvironment{corollary}{\begin{itcorollary}\rm}{\end{itcorollary}}
\newenvironment{proposition}{\begin{itproposition}\rm}{\end{itproposition}}
\newenvironment{definition}{\begin{itdefinition}\rm}{\end{itdefinition}}
\newenvironment{example}{\begin{itexample}\rm}{\end{itexample}}

\newcommand{\bcl}[1]{\begin{claim}\label{#1}}
\newcommand{\bl}[1]{\begin{lemma}\label{#1}}
\newcommand{\br}[1]{\begin{remark}\label{#1}}
\newcommand{\bt}[1]{\begin{theorem}\label{#1}}
\newcommand{\bd}[1]{\begin{definition}\label{#1}}
\newcommand{\bp}[1]{\begin{proposition}\label{#1}}
\newcommand{\bc}[1]{\begin{corollary}\label{#1}}
\newcommand{\bfact}[1]{\begin{fact}\label{#1}}
\newcommand{\bex}[1]{\begin{example}\label{#1}}
\newcommand{\bem}[1]{\begin{example}\label{#1}}
\newcommand{\ec}{\end{corollary}}
\newcommand{\eex}{\end{example}}
\newcommand{\eem}{\end{example}}
\newcommand{\el}{\end{lemma}}

\newcommand{\er}{\end{remark}}
\newcommand{\et}{\end{theorem}}
\newcommand{\ed}{\end{definition}}
\newcommand{\ep}{\end{proposition}}
\newcommand{\epr}{\end{proof}}
\newcommand{\bpr}{\begin{proof}}
\newcommand{\ecl}{\end{claim}}

\newcommand{\beq}{\begin{eqnarray}}
\newcommand{\eeq}{\end{eqnarray}}
\newcommand{\beqn}{\begin{eqnarray*}}
\newcommand{\eeqn}{\end{eqnarray*}}
\newcommand{\bi}{\begin{itemize}}
\newcommand{\ei}{\end{itemize}}
\newcommand{\ben}{\begin{enumerate}}
\newcommand{\een}{\end{enumerate}}

\newcommand{\R}{{\mathbb R}}  
\newcommand{\N}{{\mathbb N}}  

\newcommand{\url}[1]{\mathtt{#1}}

\def\Id{{1\kern-.4em 1}}

\def\epi{{\text{epi}\;}}

\def\C{{\cal C}}
\title{On the Strong Invariance Property for Non-Lipschitz
Dynamics\footnote{Corresponding Author:  Peter Wolenski. Date:
March 17, 2005.}}
\author{Mikhail Krastanov\footnote{Department of Biomathematics;
                 Institute of Mathematics and Informatics;
                 Bulgarian Academy of Sciences;
                 Acad. G. Bonchev str., bl 8;
                 1113 Sofia;
                 BULGARIA; krast@math.bas.bg }
\and  Michael Malisoff\footnote{Department of Mathematics; 304
Lockett Hall; Louisiana State University and A. \& M. College;
Baton Rouge, LA 70803-4918; USA; \{malisoff,
wolenski\}@math.lsu.edu} \and Peter Wolenski$^\ddagger$}\date{}

\begin{document}
\maketitle
\begin{abstract}
We provide a new sufficient condition for strong invariance for
 differential inclusions, under very general
conditions on the dynamics, in terms of a Hamiltonian inequality.
In lieu of the usual Lipschitzness assumption on the
multifunction, we assume a feedback realization condition that can
in particular be satisfied for measurable dynamics that are
neither upper nor lower semicontinuous.
\end{abstract}

\section{Introduction}
\label{intro} Topics in flow invariance theory provide the
foundation for considerable current research in control theory and
optimization (cf. \cite{CLR97, CLSW98, CS03, DRW04a, RW03, V00,
WZ98}).
The setting is that a multifunction $F:\R^n\rightrightarrows \R^n$
defining dynamics and a closed set $S\subseteq \R^n$ determining
state constraints are given, and the theory then contains
sufficient conditions under which some, or all, of the
trajectories of $F$ that start in $S$ remain in $S$. More
particularly, we say that $(F,S)$ is {\em weakly invariant in}
$\R^n$ provided for each $\bar x \in S$, {\em there exists} a
trajectory $t\mapsto \phi(t)$ of $F$ starting at $\bar x$ that
remains in $S$ for all $t$ before its escape time ${\rm Esc}
(\phi, \R^n)$ (precise definitions are given in section \ref{two}
below).  Weak invariance of a set is also termed {\em viability},
and sufficient (and necessary) conditions for weak invariance have
been developed under very general conditions on the dynamics, both
in terms of tangential-type inclusions and Hamiltonian
inequalities.

A more restrictive invariance property is as follows:  We say that
$(F,S)$ is {\em strongly invariant in} $\R^n$ provided for each
$\bar x \in S$, {\em each trajectory} $t\mapsto \phi(t)$ of $F$
starting at $\bar x$ remains in $S$ for all $t$ before its escape
time ${\rm Esc} (\phi, \R^n)$.
In contrast with weak invariance theorems, which merely require
the dynamics to have locally bounded values and closed graph,
sufficient conditions for strong invariance usually invoke a
Lipschitz condition on the dynamics (cf. section \ref{compare} for
a survey of results in this direction). For example, if $F$ is
locally Lipschitz and nonempty and compact-convex valued with
linear growth, then it is well known (cf. \cite[Chapter
4]{CLSW98}) that $(F,S)$ is strongly invariant in $\R^n$ if and
only if $F(x)\subseteq T^C_S(x)$ for all $x\in S$, where $T^C_S$
denotes the Clarke tangent cone.

However, this cone characterization can fail if $F$ is
non-Lipschitz, as illustrated in the following simple example:
Take $n=1$, $S=\{0\}$, $F(0)=[-1,+1]$, and $F(x)=\{-{\rm
sign}(x)\}$ for $x\ne 0$.
 Then $T^C_S(0)=\{0\}$, even though $(F,S)$ is strongly invariant in
$\R$.  This example satisfies our dynamic assumptions (cf. Example
\ref{ex0.5} below).  It is also covered by our main theorem (see
section \ref{three}).

Starting from strong invariance and its Hamiltonian
characterizations, one can develop uniqueness results and
regularity theory for proximal solutions of
Hamilton-Jacobi-Bellman equations, stability theory, infinitesimal
characterizations of monotonicity, and many other applications
(cf. \cite{AS03, CLSW98, CS03, FPR95, V00}). On the other hand, it
is well appreciated that many important dynamics are
non-Lipschitz, and may even be discontinuous, and therefore are
beyond the scope of the usual strong invariance characterizations.
Therefore, the development of conditions guaranteeing strong
invariance under less restrictive assumptions is a problem that is
of considerable ongoing research interest.

This motivates the search for sufficient conditions for strong
invariance for non-Lipschitz differential inclusions, which is the
focus of this  note.  Donchev, Rios and Wolenski \cite{DRW04a,
DRW04b} recently developed necessary and sufficient conditions for
strong invariance for so-called {\it one-sided Lipschitz}
differential inclusions.
See also \cite{RW03} for an autonomous normal type
characterization of strong invariance for certain systems with a
discontinuous component. These works apply under special
conditions on the structure of the dynamics (cf. section
\ref{three} for further details).

In this note, we pursue a very different approach. Rather than
restricting the structure of the dynamics, we provide a sufficient
condition for strong invariance under an appropriate feedback
realization hypothesis.  This hypothesis is related to Sussmann's
`unique limiting' property that was introduced in \cite{S03} in
the context of exit time optimal control problems with  continuous
dynamics, and to Malisoff's ``Lipschitz upper envelope'' condition
from \cite{M99a, M01}. Roughly speaking, our realization property
states that each trajectory $\phi$ of the dynamics $F$ is also a
unique trajectory of a nonautonomous singleton-valued dynamics $f$
for which $f(t,x)\in {\rm cone}\, \{F(x)\}$ for all $t$ and all
$x$ near $\phi(0)$ (cf. section \ref{two} below for a precise
statement of our hypothesis). This is a less restrictive
assumption than those of the known strong invariance
characterizations because it can be satisfied by important classes
of differential inclusions with measurable, but possibly neither
upper nor lower semicontinuous, right-hand sides (cf. section
\ref{two} for examples). While our main theorem can be shown
 using  Zorn's Lemma, the proof we give below is constructive, and in
 particular leads to a
 new approach to building viable
trajectories for Carath\'eodory dynamics; see Remark \ref{franko}.

In section \ref{two}, we state our realization hypothesis
precisely and provide the necessary background on differential
inclusions and nonsmooth analysis.  We also illustrate the
applicability of our hypothesis to a broad class of discontinuous
dynamics that are beyond the scope of the well known strong
invariance results.  In section \ref{three}, we announce our
strong invariance result and discuss its relationship to the known
theorems in strong invariant system theory.  Section \ref{four}
contains the proof of our strong invariance criterion, and we
close in section \ref{five} by proving a new necessary and
sufficient Hamiltonian condition for strong invariance for general
lower semicontinuous feedback realizable dynamics.

\section{Assumptions and Preliminaries}
\label{two}
\subsection{Basic Hypothesis}
Our main object of study in this note is an autonomous
differential inclusion $\dot x\in F(x)$.
 In this subsection, we state our hypothesis on $F$
and illustrate its relevance using several applications.  Our
novel feature is the requirement that each trajectory of $F$ be
realizable as the unique solution to a nonautonomous local
feedback selection of $F$. On the other hand, we will not require
the Lipschitz property or other structural assumptions on $F$ that
are generally invoked in strong invariant system theory (cf.
\cite{ CLSW98, D02, DRW04a,  FPR95, RW03}).

To make our realization hypothesis precise,  we require the
following definitions and notation. By a {\em trajectory} of $\dot
x\in F(x)$ on an interval $[0,T]$ starting at a point $x_o\in
\R^n$, we mean an absolutely continuous function $\phi:[0,T]\to
\R^n$ for which $\phi(0)=x_o$ and $\dot\phi(t)\in F(\phi(t))$ for
(Lebesgue) almost all (a.a.) $t\in [0,T]$.  We let ${\rm
Traj}_T(F,x)$ denote the set of all trajectories $\phi:[0,T]\to
\R^n$ for $F$ starting at $x$ on all possible intervals $[0,T]$,
and we set ${\rm Traj}(F,x):=\cup_{T \ge 0}{\rm Traj}_T(F,x)$ and
${\rm Traj}(F):=\cup_{x\in \R^n} {\rm Traj}(F,x)$.

A multifunction $G:\R^n\rightrightarrows\R^n$ is said to have {\em
linear growth} provided there exist positive constants $c_1$ and
$c_2$ such that $||v||\le c_1+c_2||x||$ for all $v\in G(x)$ and
$x\in \R^n$, where $||\cdot||$ denotes the Euclidean norm.  For
any interval $I$, a function $f:I\times \R^n\to\R^n$ is said to
have {\em linear growth (on $I$)} provided $x\mapsto
G(x):=\{f(t,x): t\in I\}$ has linear growth. For any sets
$D,M\subseteq \R^n$ and $\eta\in \R$, we set $M+\eta D:=\{m+\eta
d: m\in M, d\in D\}$, and ${\rm cone}\, \{D\}:=\cup\{\eta D:
\eta\ge 0\}$.  Also, ${\cal B}_n(p):=\{x\in \R^n: \Vert x-p\Vert
\le 1\}$ for all $p\in \R^n$ and ${\cal B}_n:={\cal B}_n(0)$.
 A
function $\omega(\cdot):[0,\infty)\to[0,\infty)$ is called a {\em
modulus} provided it is nondecreasing and continuous with
$\omega(0)=0$. For each $T\ge 0$,  we let $\C[0,T]$ denote the set
of all functions $f:[0,T]\times \R^{n}\to\R^{n}$ that satisfy

\begin{enumerate}
\item[$(C_1)$]  For each $x\in\R^{n}$, the map $t\mapsto f(t,x)$
is measurable;  \item[$(C_2)$]  For each compact set
$K\subseteq\R^{n}$, there exists a modulus $\omega_{f,K}(\cdot)$
such that, for all $t\in[0,T]$ and $x_1,\,x_2\in K$, $\Vert
f(t,x_1)-f(t,x_2)\Vert \le \omega_{f,K}(\Vert x_1-x_2\Vert )$; and
 \item[$(C_3)$] $f$ has linear growth on [0,T].
\end{enumerate}
 It is noteworthy that $\omega_{f,K}(\cdot)$ in the previous definition
 is
independent of $t\in[0,T]$.  For each $\bar x\in \R^n$, denote by
$\C_F([0,T], \bar x)$ those $f\in\C[0,T]$ that are also selections
of the cone of $F$ for almost all $t\in [0,T]$ and all $x\in \R^n$
sufficiently near $\bar x$; that is,
\[\C_F([0,T], \bar x):=\{f\in \C[0,T]: \exists \gamma>0 {\rm\
s.t.\ } f(t,x)\in {\rm cone}\, \{F(x)\} \text{ for a.a. }t\in[0,T]
\; \text{and all} \; x\in \gamma {\cal B}_n(\bar x)\}.\] Notice
that while elements $f\in \C_F([0,T], \bar x)$ are defined on all
of $[0,T]\times \R^n$, they are only required to  satisfy
$f(t,x)\in {\rm cone}\, \{F(x)\}$ on {\em part of} their domain.
We also let $\C_F[0,T]$ denote those $f\in \C[0,T]$ such that
$f(t,x)\in {\rm cone}\, \{F(x)\}$ for almost all $t\in [0,T]$ and
all $x\in \R^n$.  We will assume the following:
\begin{itemize}
\item[$(U)$] For each $\bar x\in \R^n$, $T\ge 0$, and $\phi\in
{\rm Traj}_T(F,\bar x)$, there exists $f\in \C_F([0,T], \bar x)$
for which $\phi$ is the unique solution   of the initial value
problem $\dot y(t)=f(t,y(t))$, $y(0)=\bar x$ on $[0,T]$.
\end{itemize}

Notice that hypothesis $(U)$ is  weaker  than requiring a
continuous selection from the dynamics $F$ that realizes the
trajectory. This is because $f$ is allowed to depend on time as
well as the state, and need only be a {\em local} selection.
Moreover, $f$ is allowed to depend on the choice of the trajectory
$\phi$, and need not be continuous. In practice, hypothesis $(U)$
can be checked using open or closed loop controls, and may be
satisfied for non-Lipschitz dynamics. The following examples
 illustrate these points and also show how to use  cones to
 check
 condition $(U)$.

\begin{example}
\label{ex0} Assume $F:\R^n\rightrightarrows\R^n$ is Lipschitz and
nonempty and compact-convex valued. We claim that $F$ satisfies
condition $(U)$. To see why, let $\bar x\in \R^n$, $T>0$, and
$\phi\in {\rm Traj}_T(F,\bar x)$ be given, and set
\[f(t,x)={\rm proj}_{F(x)}(\dot \phi(t))\] (i.e., $f(t,x)$ is
the closest point to $\dot\phi(t)$ in $F(x)$, which is well
defined by the convexity of $F(x)$).
 Then $f\in
\C_F[0,T]$ satisfies the requirement.  If on the other hand
$F:\R\rightrightarrows \R$ is  defined by $F(x)=\{1\}$ for $x<0$,
$F(0)=\{0\}\cup [1,2]$, and $F(x)=[0,2]$ for $x>0$, and if
$\phi\in {\rm Traj}(F)$, then $f(t,x)\equiv \dot \phi(t)\in {\rm
cone}\{F(x)\}$ for almost all $t$ and all $x\in \R^n$. Therefore,
condition $(U)$ is again satisfied, even though $F$ is neither
upper nor lower semicontinuous nor convex valued.
\end{example}
\begin{example}
\label{ex0.5} Consider the example from the introduction, namely,
$n=1$, $F(0)=[-1,+1]$, and $F(x)=\{-{\rm sign}(x)\}$ for $x\ne 0$.
We claim that $(U)$ is again satisfied.  To see why, let $T>0$,
$\bar x\in \R$, and $\phi\in {\rm Traj}_T(F,\bar x)$ be given.
Note that $(F,\{0\})$ is strongly invariant in $\R$.  Therefore,
either (i) $\phi$ starts at some $\bar x\ne 0$ and then moves to
$0$ at unit speed and then stays at $0$ or (ii) $\phi\equiv 0$. If
$\bar x\ne 0$, then the requirement is satisfied using
\[f(t,x)\equiv -{\rm sign}(\bar x)\beta(t),\] where $\beta(t)=1$ if
$t\in [0,|\bar x|]$ and $0$ otherwise.  In this case, $f(t,x)\in
 {\rm cone}\{F(x)\}$ for all $t\in [0,T]$ and
$x\in (|\bar x|/2){\cal B}_1(\bar x)$.   If instead $\bar x =0$,
then the requirement is satisfied with $f(t,x)\equiv 0\in {\rm
cone}\{F(x)\}$ for all $t\in [0,T]$ and $x\in \R$.
\end{example}
\begin{example}
\label{ex1} Assume $F(x)=g(x,A)U(x)$ where $A\subseteq\R^m$ is
compact, and $g:\R^n\times A\to \R^n$ is continuous in $x\in \R^n$
and measurable and satisfies
\begin{itemize}
\item[$(H)$] For each compact set $K\subseteq \R^n$, there exists
$L_K>0$ such that for all $x_1,x_2\in K$ and $a\in A$,
$(g(x_1,a)-g(x_2,a))\cdot (x_1-x_2)\le L_K\Vert x_1-x_2\Vert^2$.
Also, $x\mapsto g(x,A)$ has linear growth.
\end{itemize}
and $U:\R^n\rightrightarrows \R$ is bounded, measurable, closed
and nonempty valued, and satisfies $U(x)\cap (0,\infty)\ne
\emptyset$ and $U(x)\cap (-\infty,0)\ne \emptyset$ for all $x\in
\R^n$.
(The argument we are about to give still applies if instead of
assuming that $U(x)\cap (0,\infty)\ne \emptyset$ and $U(x)\cap
(-\infty,0)\ne \emptyset$ for all $x\in \R^n$, we instead assume
that either $U:\R^n\rightrightarrows(0,\infty)$ or
$U:\R^n\rightrightarrows(-\infty,0$).) One can easily check (cf.
\cite{BCD97}) that condition $(H)$ guarantees existence of at most
one solution $\phi:[0,T]\to \R^n$ of $\dot x=g(x,\alpha)\beta$ for
each initial condition, $T>0$, and bounded measurable functions
$\alpha:[0,T]\to A$ and
$\beta:[0,T]\to \R$. 
To check condition $(U)$, let $\phi\in {\rm Traj}(F)$. Consider
the multifunction \[G(t):=\{(a,b)\in A\times U(\phi(t)): \dot
\phi(t)=g(\phi(t),a)b\}.\]  Notice that $t\mapsto A\times
U(\phi(t))$ is closed valued and measurable, and that
$(t,(a,b))\mapsto g(\phi(t),a)b$ is measurable. Applying the
(generalized) Filippov lemma (cf. \cite[p. 72]{V00}) to $G$, we
find a measurable pair $(\alpha,\beta)$ for which $\alpha(t)\in
A$, $\beta(t)\in U(\phi(t))$, and $\dot
\phi(t)=g(\phi(t),\alpha(t))\beta(t)$ for almost all $t$. We now
show that   condition $(U)$ holds with the choice
\[f(t,x):=g(x,\alpha(t))\beta(t).\] In general, we will not have
$\beta(t)\in U(x)$ for all $t$ and $x$.  In fact, it could be that
$U(\phi(t))\cap U(x)=\emptyset$ for some $t$ and $x$, so we may
not have $f(t,x)\in F(x)$ for a.a. $t$ and $x$. On the other hand,
one can easily check that $\beta(t)\in {\rm cone}\, \{U(x)\}$ for
all $t$ and $x$, so $f(t,x)\in g(x,A){\rm cone}\{U(x)\}={\rm
cone}\{F(x)\}$ for a.a. $t$ and all $x$,
 and this gives the desired
result.
\end{example}
\begin{example}
\label{ex2} Let $R\subseteq \R^n$ be compact and
$F(x)=\{\lambda(x)+U(x)\delta(x)\}R$, where (i) $\delta,
\lambda:\R^n\to \R$ are continuous functions such that
$|\lambda(x)|\ge 2|\delta(x)|$ for all $x$, (ii)
$g(x,A):=\lambda(x)+A\delta(x)$ satisfies condition $(H)$ above
with $n=1$ and $A=[-1,+1]$, and (iii) $U:\R^n\to [-1,+1]$ is
measurable. To check condition $(U)$, let $\phi\in {\rm Traj}(F)$.
Applying the Filippov lemma as in the previous example, we find
measurable functions $u$ and $r$ such that $\dot
\phi(t)=[\lambda(\phi(t))+u(t)\delta(\phi(t))]r(t)$ for almost all
$t$, with $u(t)\equiv U(\phi(t))$.  Notice that
\[\lambda(x)+U(y)\delta(x)\in (\min\{\lambda(x)/2, 3\lambda(x)/2\},
\max\{\lambda(x)/2, 3\lambda(x)/2\})\] for all $x$ and $y$, by our
stated conditions. This implies that $\lambda(x)+U(y)\delta(x)\in
{\rm cone}\, \{\lambda(x)+U(x)\delta(x)\}$ for all $x,y\in \R^n$,
so we can satisfy condition $(U)$ using
$f(t,x):=[\lambda(x)+u(t)\delta(x)]r(t)$.  In this case $f\in
\C_F[0,T]$.
\end{example}

\subsection{Preliminaries in Nonsmooth Analysis}

The principal nonsmooth objects used in this note  are the
proximal subgradient and normal cone, and here we review  these
concepts; see \cite{CLSW98} for a complete treatment. Let
$S\subseteq \R^n$ be closed and $s\in S$.  A vector $\zeta\in
\R^n$ is called a {\em proximal normal} vector of $S$ at $s$
provided there exists $\sigma=\sigma(\zeta, s)>0$ so that
\begin{equation}\label{ProxNormal}
\langle \zeta,s'-s\rangle \leq \sigma||s'-s||^2\quad \text{for all
} s'\in S.
\end{equation}
The set of all proximal normals of $S$ at $s$ is denoted by
$N^P_S(s)$ and is a convex cone. One can show (cf. \cite[p.
25]{CLSW98}) that for each $\delta>0$ and $s\in S$, $\zeta\in
N^P_S(s)$ if and only if there exists $\sigma=\sigma(\zeta, s)>0$
so that
\begin{equation}\label{ProxNormala}
\langle \zeta,s'-s\rangle \leq \sigma||s'-s||^2\quad \text{for all
} s'\in S\cap\delta{\cal B}_n(s).
\end{equation}
Recall that the {\em distance function} $d_S(\cdot):\R^n\to\R$ is
given by $d_S(x):=\min\{\Vert x-s\Vert :s\in S\}$.

For the related functional concept, assume $f:\R^n\to
(-\infty,\infty]$ is lower semicontinuous and let $x\in{\rm
domain} (f):=\{x':f(x')<\infty\}$. Then $\zeta\in\R^n$ is called a
{\em proximal subgradient} for $f$ at $x$ provided there exist
$\sigma>0$ and $\eta>0$ so that
\begin{equation}\label{ProxSub}
f(x')\geq f(x) + \langle \zeta,x'-x\rangle - \sigma\Vert x'-x\Vert
^2\quad\text{for all  } x'\in \eta{\cal B}_n(x).
\end{equation}
The set of all proximal subgradients for $f$ at $x$ is denoted by
$\partial_P f(x)$. This set could be empty at some points, even
for $C^1$ functions (e.g., $\partial_P f(0)=\emptyset$ if
$f(x)=-|x|^{3/2}$).

We next state the version of the Clarke-Ledyaev Mean Value
Inequality  needed for our strong invariance results. Let $[x,Y]$
denote the closed convex hull of $x\in \R^n$ and $Y\subseteq
\R^n$.

\begin{theorem}\label{MVI}
Assume $x\in \R^n$, $Y\subseteq \R^n$ is compact and convex, and
$\Psi:\R^n\to(\infty,+\infty]$ is lower semicontinuous.  Then for
any $\delta<\min_{y\in Y}\Psi(y)-\Psi(x)$ and $\lambda>0$, there
exist $z\in [x,Y]+\lambda {\cal B}_n$  and $\zeta\in
\partial_P \Psi(z)$ so that
$\delta<\langle \zeta,y-x\rangle$ for all $y\in Y$.
\end{theorem}
\noindent For the proof, see  \cite[p. 117]{CLSW98}; an infinite
dimensional version also holds (see \cite{CL94}), but is not
needed here.

\subsection{Background in Differential Inclusions}
In this subsection, we review  invariant systems theory and a
standard result on compactness of trajectories for discontinuous
dynamics.  The following definition of escape times was introduced
in \cite{WZ98}:

\bd{defn2.2.1} Let $\G\in {\R}^{n}$ be open, $x_0\in \G$, and
$x(\cdot)$ be a trajectory of  a differential inclusion $\dot x\in
F(x)$ with $x(0)=x_0$ defined on a half--open interval $[0,T)$,
where $0<T\leq\infty$. Then $T$ is called an {\em escape time} of
$x(\cdot)$ from $\G$ provided at least one of the following
conditions hold:
\begin{enumerate}\addtolength{\itemsep}{-0.5\baselineskip}
\item[($E_1$)]  $T=\infty$ and $x(t)\in \G$ for all $t\geq 0$;
\item[($E_2$)]  $x(t)\in \G$ for all $t\in [0,T)$ and $\Vert
x(t)\Vert \to\infty $ as $t\uparrow T$; or \item[($E_3$)]
$T<\infty$, $x(t)\in \G$ for all $t\in [0,T)$, and
$d_{\G^c}(x(t))\to 0$ as $t\uparrow T$.
\end{enumerate}
\ed

We next define strong and weak invariance.   Assume $\G\subseteq
{\mathbb R}^{n} $ is open and $x_0\in \G$. The set of all
trajectories of $F$ originating from $x_0$ that remain in $\G$
over a maximal interval is  denoted by $\Upsilon_{(F,\G)}(x_0)$.
Therefore, $\Upsilon_{(F,\G)}(x_0)$ consists of those absolutely
continuous functions $x(\cdot)$ defined on a half--open interval
$[0,T)$ for which $x(0)=x_0$ and $\dot x(t)\in F(x(t))$ is
satisfied for almost all $t\in [0,T]$, where $T={\rm Esc}\,
(x(\cdot);\G)$.

\bd{DefInv} Assume $E \subseteq{\mathbb R}^{n}$ is closed, and
$\G\subseteq{\mathbb R}^{n}$ is open with $E\cap \G\not=
\emptyset$. Let $F:\R^n\rightrightarrows\R^n$.
\begin{enumerate}\addtolength{\itemsep}{-0.5\baselineskip}
\item[(a)]  $(F,E)$ is called {\em weakly invariant in $\G$}
provided that for every $x_0\in E\cap \G$, there exists a
trajectory $x(\cdot)\in\Upsilon_{(F,\G)}(x_0)$ that satisfies
$x(t)\in E$ for all $t\in [0,{\rm Esc}\, (x(\cdot);\G))$.
\item[(b)]  $(F,E)$ is called {\em strongly invariant in $\G$}
provided for every $x_0\in E$, every trajectory $x(\cdot)\in
\Upsilon_{(F,\G)}(x_0)$ satisfies $x(t)\in E$ for all $t\in
[0,{\rm Esc}\, (x(\cdot);\G))$.
\end{enumerate}
\ed

For Hamiltonian characterizations of strong invariance for {\em
locally Lipschitz dynamics}, see \cite{CLSW98}.  See also
\cite{DRW04a} for a  characterization of strong invariance for
systems satisfying appropriate one-sided Lipschitzness and
dissipativity conditions, and \cite{DRW04b} for general one-sided
Lipschitz dynamics (with a modified Hamiltonian).
 Our main contribution  will be a new sufficient
condition for strong invariance for dynamics satisfying the
realizability condition $(U)$, including cases where $F$ is
neither lower nor upper semicontinuous and not tractable by the
known strong invariance results.  Our condition is a Hamiltonian
inequality involving a lower semicontinuous verification function.
However,  necessity is not true generally, and it is not clear
what or if a modification of the inequality can be made to ensure
a complete characterization.

The following is a variant of the well known \lq\lq  compactness
of trajectories\rq\rq\ lemma.  This result says more than just
that a bounded set of solutions is relatively compact. Rather, a
stronger conclusion holds in that approximate trajectories have
subsequences that converge to a trajectory.  The proof is a
special case of the compactness of trajectories proof in
\cite{CLSW98}.

\begin{lemma}\label{prop2.2.2}Let $\bar
x\in \R^n$, $T>0$, $\tilde f\in \C[0,T]$ be also continuous in
$t$, and $\{y_i:[0,T]\to \R^n\}$ be a sequence of uniformly
bounded absolutely continuous functions satisfying $y_i(0)=\bar x$
for all $i$. Assume
\begin{equation}
\label{sc}
 \dot  y_i(t)\in \tilde f(\tau_i(t), y_i(t) + r_i(t)) +
\delta_i(t) {\cal B}_n \; \text{for\ a.a.} \; t\in [0,T]
\end{equation}
for all $i$,  where $\{\delta_i(\cdot)\}$ is a sequence of
nonnegative measurable functions that converges to $0$ in $L^2$ as
$i\to\infty$, $\{r_i(\cdot)\}$ is a sequence of measurable
functions converging uniformly to $0$ as $i\to \infty$, and
$\{\tau_i(\cdot)\}$ is a sequence of measurable functions
converging uniformly to $t$ on $[0,T]$ as $i\to \infty$.
  Then there exists a
trajectory $y$ of $\dot y=\tilde f(t,y)$, $y(0)=\bar x$ such that
a subsequence of $y_i$ converges to $y$ uniformly on $[0,T]$.
\end{lemma}

We will apply Lemma \ref{prop2.2.2} to continuous mollifications
of our feedback maps $f\in \C[0,T]$. More precisely, set
\[\eta(t)=
\left\{
\begin{array}{ll}
C\, {\rm exp}\left(\frac{1}{t^2-1}\right),& |t|<1\\
0, & |t|\ge 1 \end{array}\right.
\]
where the constant $C>0$ is chosen so that $\int_{\R}\eta(s)ds=1$.
For each $\eps>0$ and $t\in \R$, set
\[
\eta_\eps(t):=\frac{1}{\eps}\, \eta\left(\frac{t}{\eps}\right).
\]
Notice for future use that
\begin{equation}
\label{etaone} \int_{\R}\eta_\eps(t)dt=1\; \; \forall \eps>0.
\end{equation}
Define the following convolutions of $f\in \C[0,T]$ in the
$t$-variable:
\begin{equation}
\label{ftild} f_\eps(t,x):=\int_{\R} f(s,x)\eta_\eps(t-s)ds\;
\left(=\int_{\R} f(t-s,x)\eta_\eps(s)ds\right)
\end{equation}
with the convention that $f(s,x)=0$ for $s\not\in [0,T]$. Then
$f_\eps\in \C[0,T]$ and is continuous for all $\eps>0$.  (In fact,
$f_\eps$ is a $C^\infty$ function of $t$ for each $x\in\R^n$, but
we will not need this fact.  See \cite{E98, F99} for the well
known theory of convolutions and mollifiers.) We will apply Lemma
\ref{prop2.2.2} to a sequence $\tilde f:=f_{\eps(i)}$ with
$\eps(i)>0$ converging to zero.  In this case, we will use ideas
from the standard proof that  \[f_{\eps(i)}(\cdot,x)\to
f(\cdot,x)\; \; \text{in}\; \;  L_1[0,T] \; \; \forall x\in \R^n\]
to build trajectories of $f\in \C[0,T]$ that respect the state
constraint.

\begin{remark}
\label{cr} Note for later use that if $\tau_i(t)\equiv t$ in Lemma
\ref{prop2.2.2}, then the conclusions of the lemma remain true
even if the $t$-continuity hypothesis on $f\in \C[0,T]$ is omitted.
This follows from the proof of the compactness of trajectories
lemma in \cite{CLSW98}.\end{remark}

\section{Strong Invariance Theorem} \label{three}
\subsection{Statement of Theorem and Remarks}
Let $H_F:\R^n\times\R^n\to [-\infty, +\infty]$ denote the {\em
(upper) Hamiltonian} for our dynamics $F$; i.e.,
\[
\displaystyle H_F(x,p):=\sup_{v\in F(x)}\langle v,p\rangle.
\]
For any subset $D\subseteq\R^n$, we write $H_F(x,D)\le 0$ to mean
that $H_F(x,d)\le 0$ for all $d\in D$.  By definition, this
inequality holds vacuously if $D=\emptyset$, e.g., if $D$ is the
empty set of proximal subdifferentials of a function at some
point.

\begin{theorem}
\label{main} Let $F:\R^n\rightrightarrows \R^n$  satisfy $(U)$,
let $\Psi:\R^n\to (-\infty, +\infty)$ be lower semicontinuous, and
set ${\cal S}=\{x\in \R^n: \Psi(x)\le 0\}$. If there exists an
open set ${\cal U}\subseteq \R^n$ containing ${\cal S}$ for which
$H_F(x,\partial_P\Psi(x))\le 0$ for all $x\in {\cal U}$, then
$(F,{\cal S})$ is strongly invariant in $\R^n$.
\end{theorem}
{}We prove Theorem \ref{main} in section \ref{four} by
constructing appropriate Euler polygonal arcs;  see also Remark
\ref{franko} for an alternative nonconstructive  proof based on
Zorn's Lemma.  Theorem \ref{main} differs from the usual strong
invariance statements in the manner in which the set $S$ is
described, but it allows for some interplay between constraint and
data assumptions.   Note that we require the Hamiltonian
inequality in a {\em neighborhood} ${\cal U}$ of ${\cal S}$, for
the result is not true in general if the Hamiltonian condition is
placed only on ${\cal S}$, even if $\Psi$ and $F$ are smooth. For
example, take $n=1$, $\Psi(x)=x^2$, and $F(x)\equiv \{1\}$. In
this case, ${\cal S}=\{0\}$ and $H_F(0,\partial_P\Psi(0))=0$, but
$(F, {\cal S})$ is not strongly invariant.  On the other hand,
Example \ref{ex0.5} is covered by Theorem \ref{main}, once we
choose the verification function $\Psi(x)=x^2$. In this case, the
Hamiltonian condition reads $H_F(x,\Psi'(x))=-2x\, {\rm sign}(x)=
-2|x|\le 0$ for all $x\in \R$, so our sufficient condition for
strong invariance is satisfied.

Theorem \ref{main} contains the usual sufficient condition for
strong invariance for an arbitrary closed set $S\subseteq\R^n$ by
letting $\Psi$ be the characteristic function $I_S$ of $S$; that
is, $I_S(x)=0$ if $x\in S$ and is $1$ otherwise. Then
$\partial_P\Psi(x)=\{0\}$ for all $x\not\in {\rm boundary}\, (S)$,
and $\partial_P\Psi(x)=N^P_S(x)$ for all $x\in {\rm boundary}\,
(S)$. This implies the following special case of Theorem
\ref{main}:
\begin{corollary}
\label{cor1} Let $F:\R^n\rightrightarrows \R^n$  satisfy $(U)$ and
$S\subseteq \R^n$ be closed.  If $H_F(x,N^P_S(x))\le 0$ for all
$x\in {\rm boundary}\, (S)$, then $(F,S)$ is strongly invariant in
$\R^n$.
\end{corollary}
The converse of Corollary \ref{cor1} does not hold, as illustrated
by the simple example given in the introduction. This means that
the converse of Theorem \ref{main} does not hold.

\begin{remark}
\label{ur} Theorem \ref{main} remains true (by the same proof) if
its Hamiltonian inequality is replaced by: $\langle
f(t,x),p\rangle\le 0$ for all $T\ge 0$, $\bar x\in \R^n$, $f\in
\C_F([0,T],\bar x)$, $t\in [0,T]$, $x\in {\cal U}$, and $p\in
\partial_P\Psi(x)$.

\end{remark}

\subsection{Relationship to Known Strong Invariance Results}
\label{compare} Theorem \ref{main} improves on the known strong
invariance results because it does not require the usual Lipschitz
or other structural assumptions on the dynamics.  The papers
\cite{C75, CLR97, FPR95, K95} provide strong invariance results
for  locally Lipschitz dynamics (see also \cite[Chapter
4]{CLSW98}). In \cite{C75}, Clarke showed that strong invariance
of $(F,S)$ in $\R^n$ is equivalent to
\begin{equation}\label{tan}F(x)\subseteq T^C_S(x)\; \; \; \forall
x\in S
\end{equation}
where $T^C_S$ is the Clarke tangent cone (cf. \cite{CLSW98}).
Recall that $v\in T^C_S(x)$ if and only if for all sequences
$x_i\in S$ converging to $x$ and all sequences $t_i>0$ decreasing
to $0$, there exists a sequence $v_i\in \R^n$ converging to $v$
such that $x_i+t_i v_i\in S$ for all $i$. In particular, if
$S=\{0\}$, then $T^C_S(0)=\{0\}$.
 Later, Krastanov \cite{K95} gave an infinitesimal
characterization of normal-type, by showing strong invariance is
equivalent to the following: $H_F(x, N^P_S(x))\leq 0$ for all
$x\in S$.
 See \cite{B70, CLR97} for Hilbert space versions, and
\cite{FPR95, V00} for other strong invariance results for
Lipschitz dynamics and nonautonomous versions.

Donchev \cite{D02} extended these characterizations beyond the
autonomous Lipschitz case to \lq\lq almost continuous, one-sided
Lipschitz\rq\rq\ multifunctions.  Rios and Wolenski \cite{RW03}
proved an autonomous normal-type characterization that allows for
a discontinuous component.
 Donchev, Rios, and Wolenski
\cite{DRW04a} proved a necessary and sufficient condition for
strong invariance for a discontinuous nonautonomous differential
inclusion $F:\R^n\times {\cal I}\rightrightarrows \R^n$ whose
right-hand side is the  sum of an almost upper semicontinuous
dynamic $D(t,x)$ with nonempty compact convex values that is
dissipative in $x$, and an almost lower semicontinuous
multifunction $G(t,x)$ that is one-sided Lipschitz in $x$. In
terms of the nonautonomous Hamiltonians defined for any dynamics
$R$ by
\[
\displaystyle H_R(t,x,\zeta):=\sup_{v\in R(t,x)}\langle
v,\zeta\rangle,
\]
the main result of \cite{DRW04a} says:  If $S\subseteq \R^n$ is
closed, then $(D+G, S)$ is strongly invariant in $\R^n$ if and
only if there exists a subset $I\subseteq {\cal I}$ of full
measure in ${\cal I}$ such that
\[
H_G(t,x,N^P_S(x))-H_D(t,x,-N^P_S(x))\; \le\;  0\; \; \; \forall
(t,x)\in I\times S.
\]
This result applies to cases where the Clarke tangency condition
(\ref{tan}) is not satisfied, and covers the example in the
introduction.  They have gone further in \cite{DRW04b} to provide
a characterization of the general one-sided Lipschitz case, in
which the Hamiltonian is replaced by a limiting condition.

On the other hand, Theorem \ref{main} does not make any structural
assumptions on the dynamics.  Moreover, our feedback realizability
hypothesis $(U)$  can be satisfied for dynamics that are not
tractable by the well known strong invariance results.  For
instance, see the examples in section \ref{two}.

\section{Proof of Strong Invariance Theorem}
\label{four} This section is devoted to the proof of Theorem
\ref{main}.

Fix $T>0$ and $\bar x\in {\cal S}$.
 We first develop some properties that hold for all $f\in
\C_F([0,T], \bar x)$. Fixing $f\in \C_F([0,T], \bar x)$ and
$\eps>0$, and fixing $\gamma>0$ such that $f(t,x)\in {\rm
cone}\{F(x)\}$ for all $x\in \gamma{\cal B}_n(\bar x)$ and almost
all $t\in [0,T]$,  set
\begin{equation}
\label{Gdef}G_f^\eps[t,x,k]\; \; =\; \;
\overline{co}\left\{f_\eps(t,y): \Vert y-x\Vert \le
\frac{1}{k}\right\}\; \; \subseteq\; \; \R^n
\end{equation}
for each $t\in [0,T]$, $x\in \R^n$ and $k\in \N$, where $f_\eps$
is the regularization of $f$ defined by (\ref{ftild}) and
$\overline{co}$ denotes the closed convex hull.  This is well
defined because $f\in \C[0,T]$. By reducing $\gamma>0$ as
necessary, we can assume that $\gamma{\cal B}_n(\bar x)\subseteq
{\cal U}$.
 We also set
\[g^\eps_f[t,x,k]=1+\max\{\Vert p\Vert : p\in G_f^\eps[t,x,k]\}\]
for all $t\in [0,T]$, $x\in\R^n$, and $k\in \N$. Note that
\begin{equation}
\label{later} g^\eps_f[t,x,k]\le g_f[t,x,k]:=1+
c_1+c_2\left(||x||+\frac{1}{k}\right) \; \; \; \; \; \forall t\in
[0,T], x\in\R^n, k\in \N
\end{equation}
where $c_1$ and $c_2$ are the constants from the linear growth
requirement on $f$, so the sets $G^\eps_f[t,x,k]$ are compact. The
following  estimate is based on Theorem \ref{MVI}  from section
\ref{two}:

\bcl{claim1} If $x\in \frac{\gamma}{2}{\cal B}_n(\bar x)$, $t\ge
0$, $k\in \N$, and $h>0$ are such that
\begin{equation}\label{crucial}
0 < h \le \frac{1}{2k\, g_f[t,x,k]}\quad\text{and}\quad x+h
g_f[t,x,k]{\cal B}_n \subseteq \frac{2\gamma}{3}{\cal B}_n(\bar
x),
\end{equation}
then
\begin{equation}\label{crucialConclusion}
\Psi(x+hv)\le \Psi(x)+\frac{h}{k}
\end{equation}
holds for some $v\in G_f^\eps[t,x,k]$. \ecl

\bpr Suppose the contrary.  Fix $x\in \frac{\gamma}{2}{\cal
B}_n(\bar x)$, $t\ge 0$, $k\in \N$, and $h>0$ satisfying
(\ref{crucial}) but such that
\begin{equation}
\label{bad}  \Psi(x+hv)>\Psi(x)+\frac{h}{k}\; \; \; \forall v\in
G^\eps_f[t,x,k].
\end{equation}
  It follows  that
\begin{equation}
\label{ddef} \delta :=\frac{h}{2k}< \min_{y\in Y}\Psi(y)-\Psi(x),
\end{equation}
where $Y= x+hG^\eps_f[t,x,k]$.
 This is because $\Psi$ is lower semicontinuous.
 Let $\lambda\in (0,\frac{1}{2k})$ be such that
\begin{equation}
 \label{keyy}
x+hg_f[t,x,k]{\cal B}_n+\lambda{\cal B}_n \subseteq \gamma{\cal
B}_n(\bar x).\end{equation}
 Next we apply Theorem
\ref{MVI} with the choices $Y= x+hG^\eps_f[t,x,k]$ and $\delta$
defined by (\ref{ddef}). It follows that there exist $z \in [x,Y]
+ \lambda{\cal B}_{n}$ and $\zeta\in
\partial_P\Psi(z)$ for which
\begin{equation}\label{this}
\delta<\min_{y\in Y}\langle \zeta, y-x\rangle= \min_{v\in
G^\eps_f[t,x,k]}\langle \zeta,hv\rangle.
\end{equation}
Note that $z\in \gamma{\cal B}_n(\bar x)\subseteq {\cal U}$, by
(\ref{keyy}).
 Since $z\in [x,Y]+\lambda {\cal
B}_{n}$, (\ref{crucial}) combined with the choice of $\lambda$
gives \[\Vert z-x\Vert \le hg_f[t,x,k]+\lambda\le \frac{1}{k}.\]
Therefore $f_\eps(t,z)\in G_f^\eps[t,x,k]$, by the definition
(\ref{Gdef}) of $G^\eps_f[t,x,k]$. Since $f(s,z)\in {\rm
cone}\{F(z)\}$ for a.a. $s\in [0,T]$ (by our choice of
$\gamma>0$), our Hamiltonian hypothesis gives $\langle \zeta,
f(s,z)\rangle \le 0$ for almost all $s\in [0,T]$. Therefore,
 (\ref{this}) gives
\begin{equation}
\label{contra} 0\; \; <\; \; \delta\; \; \leq\; \; h\langle \zeta,
f_\eps(t,z)\rangle\; \; =\; \; h \int_{\R}\eta_\eps(t-s)\langle
\zeta, f(s,z)\rangle ds\; \; \le\; \; 0.
\end{equation}
The contradiction (\ref{contra}) concludes the proof of Claim
\ref{claim1}.\epr

Now set
\begin{equation}\label{KC}D:= \frac{\gamma}{2}{\cal B}_n(\bar
x)\subseteq {\cal U}.\end{equation} Let $\omega_{f,K}$ be a
modulus of continuity for $x\mapsto f(t,x)$ on $K:=D+{\cal B}_n$
for all $t\in [0,T]$. Such a modulus exists by condition $(C_2)$.
Then $\omega_{f,K}$ is also a modulus of continuity of $K\ni
x\mapsto f_\eps(t,x)$
for all $t\in [0,T]$ and $\eps>0$.
 The following estimate
follows from Carath\'eodory's Lemma (cf. \cite[p. 55]{RW98}):

\bcl{claim2}

Let $(t,x,k)\in [0,T]\times D\times \N$  and  $v\in
G^\eps_f[t,x,k]$. Then $\Vert v-f_\eps(t,x)\Vert \le
\omega_{f,K}(1/k)+1/k$.\ecl

\bpr Using the Carath\'eodory Lemma and the definition of
$G^\eps_f[t,x,k]$, we can write
\[
v=\Delta+\sum_{j=0}^{n} \alpha_{j} f_\eps(t,x_{j}),
\]
where
\[
\alpha_{j}\in [0,1]\; \; \forall j,\; \; \; \sum_{j=0}^{n}
\alpha_{j}=1, \; \; \; \Vert x_{j}-x\Vert \le 1/k \; \forall j,\;
\; \;  ||\Delta||\le 1/k.\] In particular, $x_j\in K$ for all $j$.
This gives
\[
\displaystyle
\begin{array}{lll}
||v-f_\eps(t,x)\Vert & \le & \left\Vert
v-\displaystyle\sum_{j=0}^{n} \alpha_{j}
f_\eps(t,x_{j})\right\Vert\; \; +\; \;
\left\Vert\displaystyle\sum_{j=0}^{n} \alpha_{j}
\{f_\eps(t,x_{j})-f_\eps(t,x)\}\right\Vert\; \;   \le \; \;
\displaystyle \frac{1}{k}+\omega_{f,K}\left(\frac{1}{k}\right),
\end{array}
\]
as desired.\epr

 Next define
 \[\delta(D):=1+c_1+c_2+c_2\sup\{\Vert v\Vert : v\in D\}.\]
 It follows from
 the estimate (\ref{later})  that
 \begin{equation}\label{keyforall}
 G^\eps_f[t,x,k]\subseteq \delta(D){\cal B}_n
 \; \; \forall \; t\in [0,T],\; \; x\in D,\; \; k\in \N
.\end{equation}  Next set
\begin{equation}
\label{sharp} \tilde
T:=\min\left\{T,\frac{\gamma}{8\delta(D)}\right\},\; \; \; \;
\text{and} \; \; \; h_k:=\frac{\gamma}{4k\delta(D)}\; \; \; \; \;
\forall k\in \N.
\end{equation}
Choose $N>2$
 such that
\begin{equation}
\label{fc}
 D+h_k\delta(D){\cal B}_n\; \; \subseteq\; \;  \frac{2\gamma}{3}
 {\cal B}_n(\bar x)\; \; \;
 \forall k\ge N.
\end{equation}
By the choices of $\gamma$ and $\delta(D)$,
\begin{equation}
\label{ssc} 0\; \; <\; \;  h_k \; \; \le\; \;
\frac{1}{2kg_f[t,x,k]}\; \; \; \; \; \forall t\in [0,T],\; x\in
D,\; k\in \N.
\end{equation}
Set $c(k)={\rm Ceiling}(\tilde T/h_k)$, i.e., $c(k)$ is the
smallest integer $\ge \tilde T/h_k$.
 For
each $k\ge N$, we then define a partition \[ \pi(k): 0=
t_{0,k}<t_{1,k}<\ldots< t_{c(k),k}=\tilde T \] by setting
$t_{i,k}=t_{i-1,k}+h_k$ for $i=1,2,\ldots, c(k)-1$. Since \[\tilde
T-h_k\le [c(k)-1]h_k\le \tilde T\] for all $k$, it follows that
$t_{c(k),k}-t_{c(k)-1,k}\le h_k$ for all $k$. We also define
sequences \begin{equation}\label{sequences} x_{0,k}, \; x_{1,k},
\; x_{2,k}, \ldots, \; x_{c(k),k}\; \in\; \R^n
\end{equation}
for $k\ge N$ as follows.  We set $x_{0,k}=\bar x$ and
$x_{1,k}=\bar x+(t_{1,k}-t_{0,k})v_{o,k}$, where $v=v_{o,k}\in
G^\eps_f[0,\bar x,k]$ satisfies the requirement from Claim
\ref{claim1} for the pair $(t_{0,k}, x_{0,k})=(0,\bar x)$. By
(\ref{keyforall}), we get
\begin{equation}
\label{numb} \Vert x_{1,k}-\bar x\Vert \; \le\;  h_k\delta(D)\; =
\; \frac{\gamma}{4k},
\end{equation}
so $x_{1,k}\in D$.  If $c(k)\ge 2$, then we set
  $x_{2,k}=
x_{1,k}+(t_{2,k}-t_{1,k})v_{1,k}$, where $v_{1,k}\in
G^\eps_f[t_{1,k},x_{1,k},k]$ satisfies the requirement from Claim
\ref{claim1} for the pair $(t_{1,k}, x_{1,k})$. Reapplying
(\ref{keyforall}) gives
\[ \Vert x_{2,k}-x_{1,k}\Vert \le h_k\delta(D)= \frac{\gamma}{4k},
\] so (\ref{numb}) gives \[\Vert x_{2,k}-\bar x\Vert \le
\frac{\gamma}{2k}.\] Therefore $x_{2,k}\in D$.
  We now repeat this process
using $x_{2,k}\in D$ instead of $x_{1,k}$. Proceeding inductively
 gives sequences \[v_{i,k}\in G_f^\eps[t_{i,k},x_{i,k},k]\]
 and
sequences (\ref{sequences}) that satisfy \[
x_{i+1,k}=x_{i,k}+(t_{i+1,k}-t_{i,k})v_{i,k}\] for $i=0,1,\ldots,
c(k)-1$.
 In this case,
\[
 \frac{c(k)\gamma}{4k}\le\left(\frac{\tilde
T}{h_k}+1\right)\frac{\gamma}{4k}\; =\; \tilde
T\delta(D)+\frac{\gamma}{4k}\le \frac{\gamma}{2}, \] by the
choices of $\tilde T$ and $k\ge 2$.
 Therefore,
\[ \Vert x_{i,k}-\bar x\Vert \; \le\; \frac{c(k)\gamma}{4k}\; \le
\; \frac{\gamma}{2}\]  for all $i$ and $k$.    It follows that the
sequences $\{x_{i,k}\}$ lie in $D$.

For each $k\ge N$,  we then choose $x_{\pi(k)}$ to be the
polygonal arc satisfying $x_{\pi(k)}(0)=\bar x$ and
\begin{equation}
\label{poly} \dot
x_{\pi(k)}(t)=f_\eps(\tau_k(t),x_{\pi(k)}(t)+r_k(t))+z_k(\tau_k(t))
\end{equation}
for all $t\in [0,\tilde T]\setminus \pi(k)$, where $\tau_k(t)$ is
the partition point  $t_{i,k}\in \pi(k)$ immediately preceding $t$
for each $t\in [0,\tilde T]$,
\begin{equation}
\label{pa} z_k(t_{i,k}):=v_{i,k}-f_\eps(t_{i,k},
x_{\pi(k)}(t_{i,k}))\; \; \forall i,k
\end{equation}
the $v_{i,k}\in G_f^\eps[t_{i,k}, x_{\pi(k)}(t_{i,k}),k]$ satisfy
the conclusions from  Claim \ref{claim1}, and
\begin{equation}
\label{pb}r_k(t):=x_{\pi(k)}(\tau_k(t))-x_{\pi(k)}(t)\end{equation}
for all $t\in [0,\tilde T]$ and $k$.
Then $x_{\pi(k)}$ is the polygonal arc connecting the points
$x_{i,k}$ for $i=0,1,2,\ldots, c(k)$. In particular,
$x_{i,k}\equiv x_{\pi(k)}(t_{i,k})$.

Since $x_{i,k}\in D$ for all $i$ and $k$, we conclude from Claim
\ref{claim2} that \[\sup\{\Vert z_k(\tau_k(t))\Vert : t\in
[0,\tilde T]\}\le\omega_{f,K}(1/k)+1/k \to 0\] as $k\to +\infty$.
Since $||\dot x_{\pi(k)}(t)||=||v_{i,k}||\le \delta(D)$ for all
$t\in (t_{i,k},t_{i+1,k})$ and all $i$ and $k$, and  since
\begin{equation}
\label{in} |\tau_k(t)-t|\le h_k\to 0\; \;  \text{as}\; \;  k\to
+\infty\; \; \forall t\in [0,\tilde T]\end{equation}
 we get
\[\sup\{\Vert r_k(t)\Vert: t\in[0,\tilde T]\}\to 0\; \;
\text{as}\; \;  k\to +\infty.\] Since (\ref{poly}) has the form
(\ref{sc}) from
 our
compactness of trajectories lemma and $f_\eps$ is continuous, we
can find a subsequence of $\{x_{\pi(k)}(\cdot)\}$ that converges
uniformly to a trajectory $y_\eps$ of \[\dot y=f_\eps(t,y),\; \;
y(0)=\bar x.\]  By possibly passing to a subsequence without
relabelling, we can assume that
\begin{equation}
\label{inn}
 x_{\pi(k)}\to y_\eps \; \; \text{uniformly\ \  on}\; \;   [0,\tilde
 T].\end{equation}
Moreover, since $x_{i+1,k}=x_{i,k}+(t_{i+1,k}-t_{i,k})v_{i,k}\in
D$ for all $i=0,1,\ldots, c(k)-1$ and $k\ge N$, conditions
(\ref{fc}) and (\ref{ssc}) along with  Claim \ref{claim1} give
\[
\begin{array}{cll}
\Psi(x_{1,k})&\le & \Psi(\bar x)+ \frac{h_k}{k}
\\[.8em]
\Psi(x_{2,k})&\le & \Psi(x_{1,k})+\frac{h_k}{k}
\\
&\vdots&\\
\Psi(x_{c(k),k})&\le & \Psi(x_{c(k)-1,k})+ \frac{h_k}{k}.
\end{array}\]
Summing these inequalities and recalling that $h_k\le \gamma$
gives
\[
\begin{array}{lll}
 \Psi(x_{i,k}) &\le & \Psi(\bar x)+\frac{c(k)h_k}{k}
 \\[.8em]
& \le &\Psi(\bar x)+\left(\frac{\tilde
T}{h_k}+1\right)\frac{h_k}{k}\\[.8em]
& \le &\Psi(\bar x)+\frac{1}{k}(\tilde T+\gamma)
\end{array}\]
for all $i$ and $k$ . Hence,
\begin{equation}
\label{hence} \Psi(x_{\pi(k)}(\tau_k(t)))\le \Psi(\bar
x)+\frac{1}{k}(\tilde T+\gamma)
\end{equation}
for all $t\in [0,\tilde T]$. It follows from (\ref{in}) and
(\ref{inn}) that
\begin{equation}
\label{hence2} x_{\pi(k)}(\tau_k(t))\to y_\eps(t)\; \; \forall
t\in [0,\tilde T]\; \; \text{as}\; \;   k\to
+\infty.\end{equation} Since $\Psi$ is lower semicontinuous, it
follows from (\ref{hence})-(\ref{hence2}) that $\Psi(y_\eps(t))\le
\Psi(\bar x)$ for all $t\in [0,\tilde T]$.

Now we consider a {\em sequence} $\{\eps(i)\}$ of positive numbers
converging to zero.  Let $y_i:=y_{\eps(i)}:[0,\tilde T]\to \R^n$
be the trajectories obtained by the preceding argument for
$\eps=\eps(i)$  for all $i\in \N$.
Note that $y_i(t)\in D$ for all $i$ and $t$, because each of the
polygonal arcs $x_{\pi(k)}$ constructed above joins points in $D$
and $D$ is closed and convex. Moreover,
\begin{equation}
\label{brack} \dot
y_i(t)=f(t,y_i(t))+\left[f_{\eps(i)}(t,y_i(t))-f(t,y_i(t))\right]
\end{equation}
for all $i$ and almost all $t\in [0,\tilde T]$. Since $||\dot
x_{\pi(k)}(t)||\le \delta(D)$ for all $k$ and a.a. $t\in [0,\tilde
T]$ for all the polygonal arcs $x_{\pi(k)}$ defined above, we get
\begin{equation}
\label{equi} ||y_i(t)-y_i(s)||\le \delta(D)(t-s) \; \;
\text{for}\; \;  0\le s\le t\le \tilde T \; \; \text{for\ \  all}
\; \; i.\end{equation} Since the $y_i$ are uniformly bounded and
equicontinuous, we can assume (possibly by passing to a
subsequence without relabelling) that there is a continuous
function $y:[0,\tilde T]\to D$
such that
\begin{equation}
\label{unif}
 y_i\to y\; \;  \text{uniformly\  on}\;  [0,\tilde T].\end{equation}
 We next show
that $y$ is a trajectory of $f$.  To this end, we prove the
following claim:

\begin{claim}
$f_{\eps(i)}(t,y_i(t))-f(t,y_i(t))\to 0$ in $L^2[0,\tilde T]$ as
$i\to \infty$.
\end{claim}
\begin{proof}
Since $y_i\to y$ uniformly on $[0,\tilde T]$ and $f$ is locally
bounded and locally uniformly  continuous in the $x$ variable, it
suffices to show that \begin{equation} \label{suffices}
 \tilde
\delta_i(t):=f_{\eps(i)}(t,y(t))-f(t,y(t))\to 0 \; \; \text{in}\;
\;  L^1[0,\tilde T]
\end{equation}
 as $i\to \infty$.  We do this by adapting a
standard mollification argument (see for example \cite[Chapter
8]{F99}) as follows.  We first extend $y$ to all of $\R$ by
defining $y(t)\equiv y(0)$ for all $t\le 0$ and $y(t)\equiv
y(\tilde T)$ for all $t\ge \tilde T$.  Recall that we also defined
$f(s,x)\equiv 0$ for all $s\not\in [0,T]$. It follows from
(\ref{equi}) and (\ref{unif}) that
\begin{equation}
\label{mulim}
\mu(i):=2\omega_{f,D}\left(\sup\left\{||y(t-\eps(i)z)-y(t)||:
|z|\le 1, t\in \R\right\}\right)\to 0
\end{equation}
as $i\to \infty$.  Moreover, for $0\le t\le \tilde T$, we can
change variables to get
\begin{eqnarray}
\tilde ||\delta_i(t)|| & \le  &
\int_{\R}||f(t-s,y(t))-f(t,y(t))||\eta_{\eps(i)}(s)ds\; \; ({\rm by}\; (\ref{etaone}))\nonumber\\
&=& \int_{-1}^1 ||f(t-\eps(i)z,y(t))-f(t,y(t))||\eta(z)dz
\nonumber\\
&\le &\int_{-1}^1
||f(t-\eps(i)z,y(t-\eps(i)z))-f(t,y(t))||\eta(z)dz+\mu(i)\nonumber\\
&=&\int_{-1}^1 ||g^{\eps(i)z}(t)-g(t)||\eta(z)dz +\mu(i)\nonumber
\end{eqnarray}
where \[g(t):=f(t,y(t)),\; \; \; g^{\eps(i)z}(t):=g(t-\eps(i)z).\]
  Notice that $g\in L^1[0,\tilde T]$, and that
\[||g^{\eps(i)z}-g||_1\le 2||g||_1\; \; \; \forall i.\] Applying
Minkowski's inequality for $L^1[0,\tilde T]$ therefore gives
\[
||\tilde \delta_i||_1\le \int_{-1}^1 ||g^{\eps(i)z}-g||_1\eta(z)dz
+ \tilde T\mu(i).
\]
 Moreover,
\[||g^{\eps(i)z}-g||_1\to 0\] as $i\to \infty$ for each $z\in
[-1,+1]$, by continuity of translation in the $L^1$ norm (see
\cite[Proposition 8.5]{F99}).  The desired convergence
(\ref{suffices}) therefore follows from (\ref{mulim}) and the
dominated convergence theorem.
\end{proof}

It therefore follows from Remark \ref{cr} that a subsequence of
$\{y_i\}$ converges to a trajectory of $f$ uniformly on $[0,\tilde
T]$.  This must be the aforementioned function $y$, as desired.
  Again
using the lower semicontinuity of $\Psi$, we therefore get
\[\Psi(y(t))\le \liminf_{i\to\infty}\Psi(y_{\eps(i)}(t))\le
\Psi(\bar x)\] for all $t\in [0,\tilde T]$.

The strong invariance asserted in the theorem is now immediate.
Indeed, let $x_o\in {\cal S}$, $T\ge 0$, and $\phi\in {\rm
Traj}_T(F,x_o)$ be given. We next show that
\begin{equation}
\label{goal} \bar t\; :=\; \sup\left\{t\ge 0: \Psi(\phi(s))\le
\Psi(x_o)\; {\rm for}\; 0\le s\le t\right\}\; =\; T,
\end{equation} which would imply that
 $\phi$ remains in ${\cal S}$ on $[0,T]$.  To
this end, note that if this supremum were some time $\bar t\in
[0,T)$, then the lower semicontinuity of $\Psi$ would give
\begin{equation}\label{dip1}\Psi(\phi(\bar t\, ))\; \; \le\; \;  \Psi(x_o).
\end{equation}
In particular, $\bar x:=\phi(\bar t)\in {\cal S}$. Let $f\in
\C_F([0,T-\bar t], \bar x)$ satisfy the requirement $(U)$ for $F$
and the trajectory \[[0,T-\bar t]\ni t\mapsto y(t):=\phi(t+\bar
t),\] and let $\gamma>0$ be such that $f(t,x)\in {\rm
cone}\{F(x)\}$ for almost all $t\in [0,T-\bar t]$ and all $x\in
\gamma {\cal B}_n(\bar x)$. By reducing $\gamma>0$ as necessary,
we can assume that $\gamma {\cal B}_n(\bar x)\subseteq {\cal U}$.

 By uniqueness of solutions of the initial value
problem \[\dot y=f(t,y),\; \;  y(0)=\phi(\bar t)\] on $[0,T-\bar
t]$, the first part of the proof applied to $f$ and the initial
value $\bar x=\phi(\bar t)\in {\cal S}$ would give $\tilde t\in
(0,T-\bar t\, )$ such that
\begin{equation}\label{dip2}
\Psi(\phi(\bar t+t\, ))-\Psi(\phi(\bar t\, ))\; \; \le\; \; 0 \;
\; \; \forall t\in [0,\tilde t].
\end{equation}
Here we use the fact that the trajectory on $[0,\tilde T]$
constructed above for $f$ starting at $\bar x$ can be extended to
$[0,T-\bar t]$, by the linear growth assumption $(C_3)$, and
therefore coincides with $y$ by our uniqueness assumption in (U).
 Since $\phi$ remains in ${\cal S}$ on $[0,\bar t]$, summing
(\ref{dip1})-(\ref{dip2}) would then contradict the definition of
the supremum $\bar t$. This establishes (\ref{goal}) and proves
the theorem.

{}\begin{remark} \label{franko} The preceding proof provides a
{\em constructive approach} to finding viable trajectories for our
Carath\'eodory feedback realizations $f$ that remain in ${\cal
S}$.  As suggested by \cite{F05}, an alternative but highly {\em
nonconstructive} proof of Theorem \ref{main} would proceed as
follows.  Let $x(t)$ be any trajectory of $F$ starting in ${\cal
S}$.  By Condition $(U)$, $x(t)$ admits a feedback realization
$f\in {\cal C}[0,T]$, and our Hamiltonian assumption gives
$\langle (f(t,x),0),q\rangle\le 0$ for  all $q\in
N^{\scriptscriptstyle P}_{{\scriptscriptstyle {\rm
epi}(\Psi)}}(x,\Psi(x))=\{(\xi,-1): \xi\in
\partial_P\Psi(x)\}$,
almost all $t\ge 0$, and all $x\in {\cal U}$, where
$\epi(\Psi):=\{(x,r): r\ge \Psi(x)\}$ is the epigraph of $\Psi$.
It is not hard to deduce from this (cf. \cite{CLSW98}) that
$(f(t,x),0)\in T^{\scriptscriptstyle
B}_{{\scriptscriptstyle\epi(\Psi)}}(x,r)$ for almost all $t\ge 0$,
all $x\in {\cal U}$, and all $(x,r)\in {\rm epi}(\Psi)$, where
$T^{\scriptscriptstyle B}$ denotes the Bouligand tangent cone.
Applying the measurable viability theorem (see \cite[Section
4]{FPR95}) to the dynamics $F(t,x)=\{(f(t,x),0)\}$ provides a
trajectory $t\mapsto (\phi(t),\Psi(x(0)))$ for $\dot x=f(t,x)$,
$\dot y=0$ starting at $(x(0),\Psi(x(0)))$ that stays in ${\rm
epi}(\Psi)$. This requires $f$ to be modified outside ${\cal U}$
in the usual way.  Hence, $\Psi(\phi(t))\le\Psi(x(0))\le 0$ for
all $t$.  By the uniqueness part of Condition $(U)$,
$\phi(t)\equiv x(t)$ so $x(t)$ stays in ${\cal S}$. Unfortunately,
the preceding alternative argument is highly nonconstructive,
since the proof of the measurable viability theorem relies on
Zorn's Lemma to construct the trajectory $\phi(t)$.
One natural question which should be considered is how our Euler
constructions from our proof could be used to build numerical
schemes for approximating viable trajectories for Carath\'eodory
dynamics. This type of result could be useful in physical
applications. This question will be addressed by the authors in
future research.

\end{remark}

\section{Strong Invariance Characterization}
\label{five} As we saw in Example \ref{ex0.5}, the Hamiltonian
condition that
 $H_F(x,N^P_S(x))\le 0$ for all $x\in {\rm boundary}(S)$ is
 not necessary for strong invariance for
$(F,S)$; there,  $(F,\{0\})$ is strongly invariant in $\R^n$, but
the Hamiltonian condition is not satisfied, and $F$ is upper
semicontinuous but not lower semicontinuous. On the other hand, if
we strengthen our assumption on $F$ to

\begin{itemize}
\item[]\begin{itemize} \item[\ \ $(U^\sharp)$] Condition $(U)$
holds; and $F$ is lower semicontinuous, and closed, convex, and
nonempty valued.
\end{itemize}\end{itemize}
then we get the following strong invariance characterization:
\bt{thm3}\label{las} Let $F:\R^n\rightrightarrows \R^n$ satisfy
$(U^\sharp)$ and $S\subseteq \R^n$ be closed.  Then $(F,S)$ is
strongly invariant in $\R^n$ if and only if $H_F(x,N^P_S(x))\le 0$
for all $x\in {\rm boundary}\, (S)$.\et
\begin{proof}
We showed  the sufficiency of the Hamiltonian condition for strong
invariance in Theorem \ref{main}, so it remains to show the
necessity. We do this by extending an argument from the appendix
of \cite{AS03} to non-Lipschitz $F$. Assume $(F,S)$ is strongly
invariant.  Fix $x\in {\rm boundary}(S)$, $v\in F(x)$, and
$\zeta\in N^P_S(x)$. Using Michael's Selection Theorem (cf.
\cite{RW98}), we can find a continuous selection $s:\R^n\to\R^n$
of $F$ for which $s(x)=v$. Choose $\sigma>0$ satisfying the
condition  for $\zeta$ to be in $N^P_S(x)$.

We can now use the local existence property to find $\bar t>0$ and
a trajectory $\phi:[0,\bar t]\to \R^n$ of the continuous dynamics
$y\mapsto s(y)$ starting at $x$, so
$\dot\phi(0)=s(\phi(0))=s(x)=v$. Since   $\phi\in {\rm Traj}_{\bar
t}(F,x)$ and $(F,S)$ is strongly invariant in $\R^n$, it follows
that $\phi(t)\in S$ for all $t\in [0,\bar t]$.
 This gives
$\langle \zeta, \phi(t)-x\rangle\le \sigma ||\phi(t)-x||^2$ for
all  $t\in [0,\bar t]$. Dividing  this inequality by $t\in (0,\bar
t]$, and letting $t\to 0$ gives \[\langle \zeta, v\rangle\le
\sigma \lim_{t\to 0}t||(\phi(t)-x)/t||^2=0.\]  Taking the supremum
over all $v\in F(x)$ and noting that $x\in {\rm boundary}(S)$ was
arbitrary gives the desired  result.
\end{proof}

Theorem \ref{las} is no longer true if the requirement that $F$ be
lower semicontinuous is dropped, as shown by Example \ref{ex0.5}.

\bigskip
{\small \noindent{\bf Acknowledgement:}\ This work has been funded
in part by the US National Academy of Sciences under the
Collaboration in Basic Science and Engineering (COBASE) Program,
supported by Contract No. INT-0002341 from the US National Science
Foundation. The contents of this publication do not necessarily
reflect the views or policies of the US National Academy of
Sciences or the US National Science Foundation, nor does mention
of trade names, commercial products or organizations imply
endorsement by the US National Academy of Sciences or the US
National Science Foundation. The first author (Krastanov) was also
supported by Swiss NSF Contract 7 IP 65642 (Program SCOPES) and
 by the Bulgarian Ministry of Science and Higher Education - National
 Fund for Science Research under contracts MM-807/98
 and MM-1104/01.  The second author (Malisoff) was also supported
 by the Louisiana Board of Regents Support Fund Grant
 LEQSF(2003-06)-RD-A-12.
}

\end{document}